\documentclass[a4paper,draft]{amsart}
\usepackage{amssymb}
\newtheorem{theorem}{Theorem}
\newtheorem{proposition}{Proposition}
\newtheorem{lemma}{Lemma}

\theoremstyle{definition}
\newtheorem{definition}{Definition}
\newtheorem{example}{Example}
\newtheorem{corollary}{Corollary}
\theoremstyle{remark}
\newtheorem{remark}{Remark}

\DeclareMathOperator{\MAT}{Mat}

\DeclareMathOperator{\End}{End}

\DeclareMathOperator{\Hom}{Hom}
\DeclareMathOperator{\rank}{rank}
\DeclareMathOperator{\Lie}{Lie}

\newcommand{\COMPLEXS}{\mathbb C}
\newcommand{\PGROUP}{\mathbb S}

\begin{document}

\title[Symmetrically
factorizable groups and set-theoretical
solutions]{Symmetrically factorizable groups and set-theoretical
solutions of the pentagon equation}

\author{R. M. Kashaev}
\address{
Steklov Mathematical Institute at St. Petersburg, Fontanka 27,
191011, St. Petersburg, Russia}

\email{kashaev@pdmi.ras.ru}
\author{N. Reshetikhin}
\address{Department of mathematics, University of California,
Berkeley, CA 94720, USA}
\email{reshetik@math.berkeley.edu}

\begin{abstract}
The notion of a  symmetrically factorizable Lie  group is  introduced. It
is shown that each symmetrically factorizable Lie  group is related to
a set-theoretical solution of  the pentagon  equation.  Each
simple Lie group (after a certain Abelian extension) is
symmetrically factorizable.
\end{abstract}
\maketitle

\section{Introduction}

Let $V$ be a vector space. We say that a linear operator
$S\in\End(V\otimes V)$ satisfies the  \emph{ pentagon equation},
if
\begin{equation}\label{5-gon}
S_{12}S_{13}S_{23}=S_{23}S_{12},\quad S_{ij}\in \End(V\otimes
V\otimes V)
\end{equation}
where $S_{ij}$ acts as $S$ on $i$-th and $j$-th factors in the
tensor product and leaves unchanged vectors in the remaining
factor. It is clear that the invertibility of $S$ is a necessary
property of solutions to the pentagon equation.

Let $M$ be a set. A \emph{set-theoretical} solution to the pentagon
equation is an invertible mapping
\begin{equation}\label{eq:s-map}
s:M\times M\to M\times M
\end{equation}
which satisfies the ``reversed'' pentagon equation
with respect to compositions of transformations of $M\times M\times M$
\begin{equation}\label{st-5-gon}
s_{23}\circ s_{13}\circ s_{12}=s_{12}\circ s_{23}
\end{equation}
with analogous meaning of subscripts as above.
For each transformation $s$ of $M\times M$ one associates its pull-back
$S\equiv s^*$ , i.e. the linear operator in the space
of complex valued functions on $M\times M$,
\begin{equation}\label{s-pb}
Sf(x,y)=f(s(x,y)) .
\end{equation}
It is clear that if $s$ is a set-theoretical solution to the
pentagon equation, its pull-back $S$  satisfies the operator pentagon
equation with appropriate definition for the tensor product of
infinite dimensional vector spaces.

The Heisenberg double of a Hopf algebra \cite{STS} is a natural
algebraic framework for operator pentagon equation as its
canonical element can be interpreted as a solution of the pentagon
equation in a certain algebraic sense, for example in the context
of $C^*$-algebras \cite{BS}. It can also be interpreted as the
pentagon equation for the associativity constraint in a semisimple
monoidal category generated by tensor powers of a single simple
object $X$ with $\Hom(X^{\otimes 2}, X)\simeq V$.

In this paper we study set-theoretical solutions to the pentagon
equation when $M$ is an algebraic manifold. In this case, the
pull-back solution to the pentagon equation acts on $C(M\times
M)$, the space of algebraic functions on $M\times M$. So, the
operator $S$ is defined when $s$ makes sense on a Zariski open
subset of $M\times M$. We introduce the notion of a symmetrically
factorizable Lie group and show that it provides us with
set-theoretical solutions to the pentagon equation. This paper is
one of the several papers in which we will study topological field
theories for simple complex Lie groups. In a follow-up paper we
will use the results presented here to construct corresponding
representations of the modular group of a punctured surface.

In Section~\ref{sec:1} we analyze set-theoretical solutions to the pentagon
equation. In Section~\ref{sec:2} we define the notion of symmetrically
factorizable groups and give basic examples.  Solutions of the pentagon
equation related to symmetrically factorizable Lie groups are described in
Section~\ref{sec:3}.

We would like to thank J. Bernstein for the interesting discussion.
This work was supported by CRDF grant RM1-2244 and by the NSF grant 
DMS-0070931.

\section{Set-theoretical solutions of the pentagon equation}
\label{sec:1}
\subsection{The mapping $\rho$}

Let $s$ be a transformation (invertible mapping) of $M\times M$. Define
binary operations $*$ and $\cdot$ as
\begin{equation}
  \label{eq:set-mapping}  s(x,y)=(x\cdot y,x* y) .
\end{equation}
First, we derive the equations to these operations
which ensure the pentagon equation for $s$.
\begin{proposition}[\cite{KS}]
The mapping $s$ is a set-theoretical solution to the pentagon equation if
and only if
 the following equations hold
\begin{gather}\label{eqdotstar1}
  (x\cdot y)\cdot z=x\cdot(y\cdot z)\\
\label{eqdotstar2}
(x* y)\cdot((x\cdot y)* z)=x* (y\cdot z)\\
\label{eqdotstar3}
(x* y)* ((x\cdot y)* z)=y* z
\end{gather}
and for any pair $x,y\in M$ there exists a unique pair $u,z\in M$
such that
\begin{equation}\label{eq:inv}
u\cdot z =x, \quad u*z=y .
\end{equation}
\end{proposition}

\begin{example}\label{ex:gr-case}
  If  the set $M$  is a  group with  respect to dot-operation,
$x\cdot y\equiv xy$,      then     the      only      solution to
the system~\eqref{eqdotstar2}---\eqref{eqdotstar3}  is given by
$x* y=y$ \cite{KS}. To prove this statement, notice that
invertibility of $s$ is equivalent to uniqueness of the pair
$u,z\in M$ in eqns~\eqref{eq:inv}. Using the group structure, we
can solve the first of these equations for $u$ and substitute the
result into the second. Thus, for any pair $x,y\in M$ there exists
a unique $z\in M$ such that
\begin{equation}\label{eq:xzzy}
xz^{-1}*z=y .
\end{equation}
Specifying $y=z^{-1}$ in eqn~\eqref{eqdotstar2} we have
\[
(x*z^{-1})(xz^{-1}*z)=x*1 ,
\]
which (when considered at $z=1$) implies that
\begin{equation}\label{eq:x11}
x*1=1
\end{equation}
and thus $xz^{-1}*z=(x*z^{-1})^{-1}$. Combining this with
uniqueness of $z$ for any $x,y$ in eqn~\eqref{eq:xzzy}, we
conclude that the mapping
\begin{equation}\label{eq:bij}
y\mapsto x*y
\end{equation} is a bijection for any $x\in M$.
Equation~\eqref{eqdotstar3} with $y=1$ reads
\[
(x*1)*(x*z)=1*z ,
\]
or taking into account eqn~\eqref{eq:x11}, we have
\[
1*(x*z)=1*z .
\]
Since the mapping $z\mapsto 1*z$ is a bijection as a special case of
eqn~\eqref{eq:bij}, we equivalently have
\[
x*z=z
\]
which is exactly the statement we wanted prove.
\end{example}
\begin{remark}
As is observed in \cite{kas}, in the case when $M$ is a group,
the pull-back of the corresponding set-theoretical solution of the pentagon
equation
\[
s(x,y)=(xy,y)
\]
is nothing else but the faithful realization
of the canonical element  of the Heisenberg double of the group
algebra of $M$.
\end{remark}

As one can see from this example, the essential part of the proof
that $x*y=y$ is the assumption that $x*y$ is defined for $x=1$. If
we assume that this operation is defined only on an open
dense subset of $M$ which does not  include $1$, then
(under reasonable assumptions about this dense subset) we can
still have solutions to the pentagon equation in some class of
functions.

In this paper we will focus on the case when the set $M$ is an
algebraic group. We will assume that the star operation is either
an algebraic mapping defined on a Zariski open subset of $M\times M$
or some birational correspondence. Then, equations
~\eqref{eqdotstar2}--\eqref{eqdotstar3} are expected to hold on a
Zariski open subset of $M\times M\times M$.

\begin{example}\label{exII}
Let $M=\COMPLEXS^*\times \COMPLEXS$ be group of upper triangular
matrices of the form
\[
\left(\begin{array}{cc} x_1& x_2\\ 0 &1\end{array}\right) ,
\]
We will also write $x=(x_1, x_2)$ for elements of this group. In
this notation the matrix multiplication is
\begin{equation}
  \label{eq:dot2} x\cdot y=(x_1y_1,x_1y_2+x_2) .
\end{equation}
 It was shown in \cite{KS} that this
operation together with the operation
\begin{equation}
  \label{eq:star2}
x* y=(y_1x_2(x_1y_2+x_2)^{-1},y_2(x_1y_2+x_2)^{-1})
\end{equation}
defined on a Zariski open subset of $M\times M$ satisfies
equations ~\eqref{eqdotstar1}--\eqref{eqdotstar3}.
\end{example}
Notice that this operation differs from $x*y=y$ and the reason is
that it is not defined for $x=1$ and, therefore, arguments used in
the previous example do not apply.

The following proposition gives a construction of the star
operation if the dot operation comes from a group with a special
structure.

\begin{proposition}[\cite{KS}]\label{prop:rho}
Let $M$ be an algebraic group. Assume that it is equipped with
invertible mapping (or birational correspondence)
\[
\rho\colon M\to M
\]
defined on a Zariski open subset of $M$ satisfying the condition that
\begin{equation}
  \label{eq:rho}
\sigma(x)\equiv(\rho(y))^{-1}\rho\left(\rho(x)(\rho(xy))^{-1}\right)
\end{equation}
is independent of $y$ on a Zariski open subset of $M\times M$.
Then, operations
\begin{equation}\label{eq:s}
x\cdot y=xy, \ x*y=\rho(x)(\rho(xy))^{-1}
\end{equation}
define  a set-theoretical solution of the pentagon equation.
\end{proposition}
\begin{proof}
  Equations (\ref{eqdotstar1})---(\ref{eqdotstar3})    can    be    checked
straightforwardly.      The    inverse     mapping to $s$ defined
by (\ref{eq:s})  has the form
\[
(x,y)\mapsto
\left(\rho^{-1}(y\rho(x)),(\rho^{-1}(y\rho(x)))^{-1}x\right).
\]
\end{proof}

\begin{corollary}
If $\rho(x)$ satisfies the conditions of
Proposition~\ref{prop:rho} then so does $\tilde\rho(x)=\rho(x)a$
for any $a\in M$ and they define one and the same solution of the
pentagon equation.
\end{corollary}
This corollary demonstrates that if the mapping $\rho$ exists it
is not unique. Fixing the equivalence class of $\rho$,  one can
choose having in mind some additional properties. In the simplest
case of Example~\ref{ex:gr-case} one can choose $\rho(x)=x^{-1}$.
In other examples we will have $\rho$ of order three and will use
it later in a construction of the representation of the mapping
class group of a surface.

\begin{example}\label{ex:2rho}In Example~\ref{exII} the order three
mapping $\rho$ defined as
\begin{equation}\label{eq:ro2}
\rho(x)=(x_2x_1^{-1},x_1^{-1})
\end{equation}
satisfies the conditions of Proposition~\ref{prop:rho} with
$\sigma(x)=\rho^{-1}(x)$.
\end{example}
\begin{example}\label{ex:3}
  Let us fix rational numbers $s_1,s_2,s_3$ such that
  \[
s_1+s_2+s_3=0,\quad s_2(s_3-s_1)\ne0
\]
  and define
\[
M\equiv\{x=(x_0,x_1,x_2,x_3)\in\COMPLEXS^4\}
\]
with dot-mapping
\[
x\cdot y=\left(x_0y_0,x_1y_0^{s_{21}}+y_1,
x_2y_0^{s_{31}}+y_2+x_1y_3y_0^{s_{21}},x_3y_0^{s_{32}}+y_3\right),\quad
s_{ij}\equiv s_i-s_j
\]
given by the
multiplication rule of triangular $3\times3$-matrices
\[
x\mapsto\left(
\begin{array}{ccc}
x_0^{s_1}&x_0^{s_1}x_1&x_0^{s_1}x_2\\ 0&x_0^{s_2}&x_0^{s_2}x_3\\ 0&0&x_0^{s_3}
\end{array} \right)
\]
The following birational correspondence
\begin{equation}\label{eq:rho-sl3}
\rho(x)=\left(
\left(x_2^{s_3}x_4^{s_1}\right)^{1/(s_2s_{31})},
  x_0^{s_{32}}x_1x_4^{-1},x_0^{s_{31}}x_4^{-1},
 x_0^{s_{21}}x_3x_2^{-1}\right)
\end{equation}
where
\(
x_4\equiv x_1x_3-x_2
\), is of order three and
satisfies eqn~(\ref{eq:rho}) with $\sigma(x)=\rho^{-1}(x)$.
\end{example}

\subsection{Set-theoretical solutions and group operations}

If $s$ is a set-theoretical solution to the
pentagon equation then so is the mapping
\(
\bar{s}=s_{21}^{-1}
\).
Using this fact, we define two more binary operations on $M$:
\begin{equation}
  \label{eq:odot-ast}   \bar   s(x,y)\equiv(x\odot y,x\circledast y)
\end{equation}
Clearly, the  pentagon equation implies that these operations also
satisfy relations~(\ref{eqdotstar1})---(\ref{eqdotstar3}) with replacements
\[
\cdot\to\odot,\quad *\to\circledast
\]
In particular, $\odot$ is an associative operation. Besides that, from
the definition of $\bar s$ it follows that
\begin{equation}
  \label{eq:s-sbar}    (x*y)\odot(x\cdot    y)=y,\quad
(x*y)\circledast(x\cdot y)=x
\end{equation}
and
\begin{equation}
\label{eq:s-bar-op} (x\circledast y)\cdot(x\odot    y)=y,\quad
(x\circledast y)*(x\odot y)=x .
\end{equation}
In fact the first set of equations is equivalent to the last one.
Notice a symmetry of the equations for these operations with respect
to replacements $(\cdot,*)\leftrightarrow(\odot,\circledast)$.

Let $G$ be an algebraic group and the operation $x\cdot y=xy$ be
the group multiplication in $G$. Denote by $i$ the operation of
taking inverse in $G$. If we are looking for a set-theoretical
solution to the pentagon equation with $M=G$, we should assume
that associative operation $\odot$ is defined on a Zariski open
subset of $G$. Assume that it is ``almost a group'', i.e. that
there exists involutive birational correspondence $j: G\to G$
defined on a Zariski open subvariety of $G$ such that neither
$x\odot j(x)$ nor $j(x)\odot x$ are necessarily defined, but the
identities
\[
(y\odot x)\odot j(x)=j(z)\odot (z\odot y)=y
\]
hold.

\begin{theorem}\label{prop:jot}
Let $G$, $s$, $\bar s$ be as above.   Mapping $s$ is a set-theoretical
solution of the pentagon equation if and only if
\begin{equation}
  \label{eq:jot2}   j(x\cdot   y)=j(x)\cdot  j(k(x)\cdot   j(y)),\quad
 k\equiv i\circ j\circ i
\end{equation}
on a Zariski open subset of $G\times G$. Furthermore,
eqn~\eqref{eq:jot2} implies that the mapping
\begin{equation}\label{eq:rho-i-j}
\rho=i\circ j
\end{equation}
defined on a Zariski open subset of $G$ is of order three and
satisfies all conditions
of Proposition~\ref{prop:rho} with $\sigma(x)=\rho^{-1}(x)$.
\end{theorem}
\begin{proof}
  To  simplify formulae we  will omit  the dot  in denoting  the group
multiplication in $(G,\cdot,i)$ as  well as the composition symbol for
mappings $i$ and $j$.

  Suppose eqn~(\ref{eq:jot2})  holds. Using  the same equation  in its
right hand side we obtain the consistency condition
\[
j(xy)=j(x)jiji(x)j(xy)
\]
which is equivalent to
\begin{equation}
  \label{eq:jot1} k\equiv iji=jij.
\end{equation}
This implies  that mapping $\rho$,  defined in eqn~(\ref{eq:rho-i-j}),
is of third order.  Expression~(\ref{eq:rho}) is evaluated as follows
\[
  i\rho(y)\rho(\rho(x)i\rho(xy))=j(y)\rho(ij(x)j(xy))
=j(y)\rho j(iji(x)j(y))=ji(x)=\rho^{-1}(x).
\]
Thus, all  conditions of Proposition~\ref{prop:rho}  are satisfied with
 $\sigma(x)=\rho^{-1}(x)$.

Assume   now  that $s$ is a set-theoretical solution of the pentagon equation.
 The first of eqns~(\ref{eq:s-sbar})
can be solved for the $*$-operation:
\[
x*y=y\odot j(xy) .
\]
Substituting it into (\ref{eqdotstar2}), we obtain
\begin{equation}
  \label{eq:dot-odot1} (y\odot j(j(x)i(z)))(z\odot x)=(yz)\odot x,
\end{equation}
where  we   have  replaced $x\to j(x)i(yz)$. Substituting here
$x\to j(x)$, then multiplying from the right by $x$ with respect
to $\odot$, and replacing consequently $y\to y\odot(xi(z))$ and
$z\to z\odot x$, we obtain
\begin{equation}
  \label{eq:odot-dot1} (y\odot(x i(z\odot x))(z\odot x)=(yz)\odot x
\end{equation}
where we have also exchanged the left and right hand sides.
Comparing eqns~(\ref{eq:dot-odot1}) and (\ref{eq:odot-dot1}) we obtain
the relation
\[
j(j(x)i(z))= x i(z\odot x)
\]
which can be solved for the operation $\odot$:
\begin{equation}
  \label{eq:more-dot} x\odot y=ij(j(y)i(x))y .
\end{equation}
The  symmetry $(\cdot,i)\leftrightarrow(\odot,j)$ leads to the
equation
\[
x y=ji(i(y)\odot j(x))\odot y ,
\]
which after using eqn~(\ref{eq:more-dot}) gives
\[
x y=ji(ij(xy)j(x))\odot y=j(ij(x)j(xy))\odot
y=ij(j(y)ij(ij(x)j(xy)))y.
\]
Solving this for $j(xy)$, we come to equation~(\ref{eq:jot2}).
\end{proof}

\begin{example}\label{ex:2} Let $M=\COMPLEXS^*\times \COMPLEXS$ with the
dot and star operations being given as in Example~\ref{exII}. The
$\odot$ operation can be found from the equation
(\ref{eq:s-sbar}):
\[
x\odot y=(x_2y_1 +x_1,x_2y_2).
\]
Clearly, this operation is equivalent to the matrix multiplication
of $2\times 2$ matrices
\[
\left(\begin{array}{cc} 1&0 \\x_1 & x_2\end{array}\right).
\]
The inverse for the multiplications $\cdot$ is
\[
i(x)=\left(1/x_1, -x_2/x_1\right).
\]
The operation $j$ exists and is given by
\[
j(x)=\left(-x_1/x_2,1/x_2\right).
\]
If we include the line $(0,y)$ in $G$ and remove the
line $(x,0)$ then this operation becomes the inverse for the
associative operation $\odot$ with the unit $(0,1)$.

Composition of $i$ and $j$ operations give
\[
k(x)=iji(x)=jij(x)=(x_2,x_1) ,
\]
while the third order mapping
\[
\rho(x)=i(j(x))=\left(-x_2/x_1, 1/x_1\right)
\]
differs from that of Example~\ref{ex:2rho} , see
eqn~\eqref{eq:ro2}, in the sign of the first component. This is
equivalent to taking the composition with the automorphism $x_i\to
-x_i$.
\end{example}

\begin{corollary}
  Eqn~(\ref{eq:jot2})
  implies that the group $(G,\cdot,i)$  is birationally equivalent
  to "almost a group" $(G,\odot,j)$. The equivalence
  is given by the mapping $k=iji$:
\[
x\odot y=k(k(x)k(y)) .
\]
\end{corollary}
\begin{proof}
 Indeed,   let    us   rewrite   eqn~(\ref{eq:more-dot})    by   using
 (\ref{eq:jot2})
\[
  x\odot                            y=i(yj(kj(y)ji(x)))y=ij(kj(y)ji(x))
=iji(iji(x)ikj(y))=k(k(x)k(y)).
\]
Intertwining    of     the    inverse    operation     follows    from
eqn~(\ref{eq:jot1}):
\(
j=iki=kik
\).
\end{proof}

\begin{corollary}
  Operations $i$ and $j$ realize the action of the
permutation group $\PGROUP_3$ on $M$.
\end{corollary}
\begin{proof}
It is  clear that $i^2=j^2=1$ which  together with eqn~(\ref{eq:jot1})
imply that  $(ij)^3=1$.  Thus, elements  $i$ and $j$ generate  a group
which is isomorphic to $\PGROUP_3$.
\end{proof}

\begin{remark}This $\PGROUP_3$  symmetry reflects the tetrahedral symmetry
   of the corresponding set-theoretical solution of the pentagon equation.
\end{remark}

\section{Symmetric factorization in a Lie group}
\label{sec:2}
Here we introduce the  notion  of  \emph{symmetric
factorization} in a Lie group.
\begin{definition}
A Lie group $G$ with the Lie algebra ${\mathfrak g}$ is called
\emph{symmetrically factorizable} if
\begin{itemize}
\item ${\mathfrak g}$
as a vector space is isomorphic to the direct sum of its Lie
subalgebras ${\mathfrak g}_+$ and ${\mathfrak g}_-$, i.e.
 $\mathfrak g=\mathfrak g_+\oplus\mathfrak g_-$.
\item there exists an element $\theta\in G$ which conjugates
Lie subgroups $G_\pm$ corresponding to Lie subalgebras ${\mathfrak
g}_\pm$, i.e. $\theta G_-=G_+\theta$.
\end{itemize}
\end{definition}

The element $\theta$ is called \emph{conjugating element}. In a
symmetrically factorizable Lie group there exists an open dense
neighborhood $G'\subset G$ of the unit element such that for any
$g\in G'$ there exist $g_\pm, \bar{g}_\pm \in G_\pm$ such that
\begin{equation}\label{factor}
g=g_+ g_-^{-1}=\bar{g}_-^{-1}\bar{g}_+ .
\end{equation}
\begin{example}\label{gl-2n}
Let $G=GL(2N,\COMPLEXS)$. Any element of this Lie group can be
regarded as a block matrix of the form
\[
g=\left(\begin{array}{cc} g_{11}&g_{12}\\ g_{21}&g_{22}
\end{array}\right),\quad \det(g)\ne 0,\quad
g_{ij}\in\MAT(N,\COMPLEXS).
\]
Subgroups
\[
G_+=\{g\in  G|\  g_{21}=0,\  g_{22}=1\},\quad  G_-=\{g\in G|\
g_{12}=0,\ g_{11}=1\}
\]
are conjugate to each other
\[
\theta    G_-=G_+\theta,\quad   \theta\equiv\left(
\begin{array}{cc} 0&b\\   b^{-1}&0   \end{array}
\right),\quad   b\in GL(N,\COMPLEXS) .
\]
Thus,  Lie group $GL(2N,\COMPLEXS)$ is
symmetrically factorizable.

Let $G'\subset G$ be the subset of all elements $g$ such that
\[
g_{22},g_{11}\in GL(N,\COMPLEXS) .
\]
Every $g\in G'$ has
\emph{unique} factorization:
\[
 g=g_+g_-^{-1}=\bar{g}_-^{-1}\bar{g}_+ ,
\]
where
\[
g_+=\left(\begin{array}{cc}
    g_{11}-g_{12}g_{22}^{-1}g_{21}&g_{12}g_{22}^{-1}\\ 0&1
\end{array}\right),\quad
g_-^{-1}=\left(\begin{array}{cc} 1&0\\ g_{21}&g_{22}
\end{array}\right)
\]
and
\[
\bar{g}_-^{-1}=\left(\begin{array}{cc} 1&0\\
g_{21}g_{11}^{-1}&g_{22}- g_{21}g_{11}^{-1}g_{12}
\end{array}\right), \quad \bar{g}_+=\left(\begin{array}{cc}
    g_{11}&g_{12}\\ 0&1
\end{array}\right) .
\]
Thus, we obtain explicit description of the symmetrical
factorization of  $GL(2N,\COMPLEXS)$.
\end{example}
\begin{example}\label{simple-alg}
 Let $G$ be  a finite dimensional complex algebraic
 simple Lie  group with fixed Borel
 subgroup $B\in G$. On  a Zariski open  subset   $G'\subset G$  we  have  the  Gauss decomposition
\[
G'= N_+HN_-\cap N_-HN_+ ,
\]
where $N_\pm$  are nilpotent in  the Borel subgroups
$B_\pm=HN_\pm$ ($B_+\equiv B$). The  Weyl group  $W$ of $G$  can be
naturally identified with the quotient group
\[
W= N(H)/H,\quad N(H)=\{g\in G\colon gHg^{-1}=H\} .
\]
Let $w_0\in W$ be the longest element in $W$ and
$\{h_i\}_{i=1}^r$, $r=\rank(G)$, be a basis in the Cartan
subalgebra ${\mathfrak h}\in {\mathfrak g}=\Lie(G)$ generated by
simple roots. We choose an enumeration of simple roots which for
root system of type $A_n, D_n, E_6$ corresponds to the Cartan
matrices
\[
\begin{array}{lc}
A_n:& \left(\begin{array}{ccccc}  2 & -1 & 0 &\dots & 0 \\
                                 -1&  2 & -1 &\dots & 0 \\
                                 0& -1 & 2 & \dots & 0 \\
                                 \vdots &\vdots&\vdots &\ddots&\vdots \\
                                 0 & 0 & 0 & \dots & 2
\end{array}\right)\\
D_n:& \left(\begin{array}{ccccccc} 2 & -1 & 0 &\dots & 0 & 0 & 0\\
                                 -1&  2 & -1 &\dots & 0 & 0& 0\\
                                 0& -1 & 2 & \dots & 0 & 0 & 0\\
                          \vdots&\vdots&\vdots &\vdots&\ddots &\vdots&\vdots \\
                                 0&0&0& \dots &2& -1 & -1\\
                                 0&0&0& \dots &-1 & 2 & 0 \\
                                 0&0&0& \dots &-1 & 0& 2
\end{array}\right)\\
E_6:& \left(\begin{array}{cccccc} 2 &-1& 0&0&0&0\\
                                 -1&2&-1&0&0&0 \\
                                 0&-1&2&-1&-1&0 \\
                                 0&0&-1&2&0&0 \\
                                 0&0&-1&0&2&-1\\
                                 0&0&0&0&-1&2
\end{array}
\right)
\end{array}
\]

The longest element of the Weyl group acts on $h_i$ as:
\[
w_0(h_i)=-h_{\tau(i)} ,
\]
where $\tau$ is an automorphism of order 2 of the set of vertices
of the Dynkin diagram of the Lie algebra $\mathfrak g$. The
automorphism $\tau$ is non-trivial only for root systems of the
type $A_n, D_n, E_6$ where it acts as follows,
\begin{equation}\label{tau}
\begin{array}{ll}
 A_n:& \tau (i)= n+1-i,\quad i=1,\dots, n ;\\
 D_n:& \tau (i)= \left\{
\begin{array}{ll}
n&\mbox{if }\  i=n-1,\\
n-1&\mbox{if }\  i=n,\\
i&\mbox{otherwise;}
\end{array}\right.\\
 E_6:& \tau (i)= \left\{
\begin{array}{ll}
6&\mbox{if }\  i=1,\\
1&\mbox{if }\  i=6,\\
i&\mbox{otherwise;}
\end{array}\right.
 \end{array}
\end{equation}

Consider the decomposition of the Cartan subalgebra $\mathfrak h$
of the Lie algebra $\mathfrak g$ into a direct sum
\[
\mathfrak g\supset {\mathfrak h}={\mathfrak h}_0\oplus{\mathfrak
h}'\oplus{\mathfrak h}'' ,
\]
where ${\mathfrak h}_0$ is the subspace on which $w_0$ acts as
$-1$, while  $w_0({\mathfrak h}')={\mathfrak h}''$. The following
table gives the dimensions of the subspaces:
\[
\begin{array}{|c|c|c|}
\hline
\mathfrak g&\dim(\mathfrak h_0)&\dim(\mathfrak h')=\dim({\mathfrak h}'')\\
\hline
\hline
X_n\ne A_n, D_n, E_6&n&0\\
\hline
A_{2k}&0&k\\
\hline
A_{2k+1}&1&k\\
\hline
D_n&n-2&1\\
\hline
E_6&2&2\\
\hline
\end{array}
\]
Denote by $H_0, H', H''$ the subgroups of the Cartan subgroup $H$
corresponding to ${\mathfrak h}_0,{\mathfrak h}',{\mathfrak h}''$,
respectively, and consider a Lie group\footnote{In the case of $A_{2k}$
we take just $D=G$.} $D=G\times H_0$ together with its
subgroups
\[
\begin{array}{lcl}
D_+&=&\{(h_0h'f, h_0)|\ h'\in H',\ h_0\in H_0,\ f\in N_+\},\\
D_-&=&\{(h_0^{-1}h''f, h_0)|\ h''\in H'',\ h_0\in H_0,\ f\in
N_-\}.
\end{array}
\]
Let $\bar{w}_0\in N(H)\subset G$ be a representative of the
longest element $w_0\in W$ of the Weyl group. Set $\theta=(\bar
w_0,1)$. We have
\[
\theta(h_0^{-1}h''f, h_0)\theta^{-1}=(h_0w_0(h'')\bar{w_0}f
\bar{w_0}^{-1}, h_0)\in D_+ .
\]
Thus, the triple $(D, D_-, D_+)$ determines a factorizable Lie
group with the conjugating element $\theta$. Notice that we can
always choose  $\bar w_0$ such that $\bar w_0^2$ is central. For
this choice $\theta^2$ is central in $D$. If necessary, taking
quotient with respect to a center of $G$, we can always assume
that $\theta^2=1$.
\end{example}

\section{Symmetrical factorization and the
pentagon equation}
\label{sec:3}
Now  we   can  relate  symmetrically  factorizable   Lie  groups
to set-theoretical solutions of the pentagon equation. Our basic
examples are symmetrically factorizable Lie groups related to
simple algebraic Lie groups. In this case we know (see examples
above) that we can assume that the  conjugating element $\theta$ can
be chosen such that
$\theta^2=1$. In fact, this choice can be made under more general
assumptions.
\begin{proposition}
If element $\theta^2$ is uniquely left factorizable,
\[
\theta^2=(\theta^2)_+(\theta^2)_-^{-1}
\]
then
\[
\tilde\theta\equiv(\theta^2)_+^{-1}\theta
\]
is also a conjugating element and it is unipotent.
\end{proposition}
\begin{proof} It is clear that $\tilde\theta$ conjugates $G_\pm$
\[
\tilde\theta G_-\tilde\theta^{-1}=(\theta^2)_+^{-1}\theta
G_-\theta^{-1}(\theta^2)_+ =(\theta^2)_+^{-1}G_+(\theta^2)_+=G_+ .
\]
The unipotency of $\tilde\theta$ follows from the equalities:
\begin{multline*}
\tilde\theta^2=(\theta^2)_+^{-1}\theta(\theta^2)_+^{-1}\theta=
(\theta^2)_+^{-1}\theta^2\theta^{-1}(\theta^2)_+^{-1}\theta=
(\theta^2)_-^{-1}\theta^{-1}(\theta^2)_+^{-1}\theta\\
=(\theta^{-1}(\theta^2)_+\theta(\theta^2)_-)^{-1}
=\left(\theta^2\theta^{-1}(\theta^2)_+^{-1}\theta\right)_-^{-1}
=\left(\theta(\theta^2)_+^{-1}\theta^2\theta^{-1}\right)_-^{-1}\\
=\left(\theta(\theta^2)_-^{-1}\theta^{-1}\right)_-^{-1}=1 .
\end{multline*}
Notice that uniqueness of factorization of $\theta^2$ is essential here.
\end{proof}
\begin{lemma}\label{le:1}
If both $\theta$ and $\theta^{-1}$ are conjugating  elements, 
then for any $g\in G$ the
following relations hold:
\begin{equation}\label{eq:g-1}
\left(\theta g^{-1}\theta^{-1}\right)_\pm=\theta g_\mp\theta^{-1}.
\end{equation}
\end{lemma}
\begin{proof}
In this case we have simultaneously $\theta G_\pm=G_\mp\theta$,
so that for any $g\in G$
\[
\theta g^{-1}\theta^{-1}=\theta g_-\theta^{-1} \left(\theta
g_+\theta^{-1}\right)^{-1} .
\]
\end{proof}
The next theorem can be regarded as the main result of the paper.
\begin{theorem}
\label{prop:g-plus-minus} (A) Let $G$ be a symmetrically
factorizable algebraic group. Then the mapping
\[
\rho(x)\equiv(\theta x)_+^{-1},\quad \rho^{-1}(x)=(\theta^{-1}x^{-1})_+
\]
defined  on  a Zariski open subset of $G_+$,   satisfies  the   conditions  of
Proposition~\ref{prop:rho}.

(B) If furthermore $\theta^2=1$, then
\(
j(x)\equiv (\theta x)_+
\)
is involutive and satisfies eqn~(\ref{eq:jot2}).

(C) With any involutive solution $j$ of  eqn~(\ref{eq:jot2}) there
exists an associated symmetrically factorized ``almost a group''.
\end{theorem}
\begin{proof}
  (A)
Evaluating (\ref{eq:rho}) we have
\begin{multline*}
(\rho(y))^{-1}\rho\left(\rho(x)(\rho(xy))^{-1}\right)=(\theta y)_+\left(
\theta(\theta x)_+^{-1}(\theta xy)_+\right)_+^{-1}
=(\theta y)_+\left(\theta(\theta x)_+^{-1}\theta xy\right)_+^{-1}\\
=(\theta y)_+\left(\theta(\theta x)_-^{-1}y\right)_+^{-1}
=(\theta y)_+\left(\theta(\theta x)_-^{-1}\theta^{-1}(\theta y)_+\right)^{-1}
=\theta(\theta x)_-\theta^{-1}
\end{multline*}
which is independent of $y$.

(B) Mapping  $j$ is  involutive:
\[
j\circ  j(x)=(\theta(\theta  x)_+)_+=(\theta^2 x)_+=(x)_+=x
\]
while  eqn~(\ref{eq:jot2}) is checked as follows
\begin{multline*}
j(xy)=(\theta xy)_+=\left((\theta x)_+(\theta x)_-^{-1}y\right)_+=
(\theta x)_+\left((\theta x)_-^{-1}\theta(\theta y)_+\right)_+
\\
=j(x)j(\theta(\theta x)_-^{-1}\theta j(y))=j(x)j(k(x)j(y))
\end{multline*}
where by using Lemma~\ref{le:1}
\[
k(x)=\theta(\theta x)_-^{-1}\theta =(\theta(\theta x)^{-1}\theta)_+^{-1}=
(\theta x^{-1})_+^{-1}=iji(x)
\]

(C)  ``Almost a group'' $G$ is $(M\times M)\cup \theta$ as a set
with the following multiplication table
\begin{gather*}
(x_1,x_2)(y_1,y_2)\equiv(x_1j(x_2j(y_1)),k(k(x_2)y_1)y_2),
\quad \theta^2=(1,1),\\
\theta(x_1,x_2)=(j(x_1),k(x_1)x_2),\quad(x_1,x_2)\theta=(x_1j(x_2),k(x_2))
\end{gather*}
where $k\equiv iji$.
The unit element and the inverse operation are correspondingly
\[
1\!\!1=(1,1),\quad (x_1,x_2)^{-1}=(ik(k(x_1)x_2),ij(x_1j(x_2)))
\]
The subgroups $G_\pm$ are
\[
G_+=\{(x,1)|\ x\in M\},\quad G_-=\{(1,x)|\ x\in M\}
\]
with the conjugating element $\theta$:
\[
\theta(1,x)=(x,1)\theta, 
\]
and the factorization formula
\[
(x_1,x_2)=(x_1,1)(1,x_2).
\]
\end{proof}
\begin{corollary}\label{cor:cen}
If $\theta^2$ is central, then $\sigma(x)=\rho^{-1}(x)$ and $\rho$
is of order three.
\end{corollary}
\begin{proof}
In this case we have simultaneous conjugations
$\theta G_\pm=G_\mp\theta$ so that
\begin{multline*}
\rho\circ\rho(x)=\left(\theta(\theta x)_+^{-1}\right)_+^{-1}=
\left(\theta(\theta x)_-^{-1}x^{-1}\theta^{-1}\right)_+^{-1}\\
=\left(\theta(\theta x)_-^{-1}\theta^{-1}\theta
x^{-1}\theta^{-1}\right)_+^{-1} =\left(\theta(\theta
x)_-^{-1}\theta^{-1}\right)^{-1} =\theta(\theta
x)_-\theta^{-1}=\sigma(x) .
\end{multline*}
Using Lemma~\ref{le:1}, we have
\[
\sigma(x)=\theta (\theta x)_-\theta^{-1}=
\left(\theta (\theta x)^{-1}\theta^{-1}\right)_+=
\left(\theta x^{-1}\theta^{-2}\right)_+=
\left(\theta^{-1} x^{-1}\right)_+=\rho^{-1}(x)
\]
Thus, $\rho^3(x)=x$.
\end{proof}

\begin{example}\label{ex:gl2} Consider the Lie group $G=GL(2,\COMPLEXS)$ with
the symmetrical factorization (see example \ref{gl-2n}) determined
by subgroups
\[
G_+=\{\left(\begin{array}{cc} x_1&x_2 \\0 &1\end{array}\right)|\
x_1\in \COMPLEXS^*,\  x_2\in \COMPLEXS\}
\]
\[
G_-=\{\left(\begin{array}{cc} 1&0 \\x_1 &x_2\end{array}\right)|\
x_1\in \COMPLEXS, \ x_2\in \COMPLEXS^*\}
\]
with the conjugating element
\[
\theta=\left(\begin{array}{cc}0&b\\c&0\end{array}\right), \quad
bc\in\COMPLEXS^* .
\]
It is easy to verify that if we parameterize $g\in G_+$ by coordinates
$x_1,x_2$ as above, then
\[
\rho(x)=(\theta
g)_+^{-1}=\left(\begin{array}{cc}-x_2/(bx_1)&1/(cx_1)\\0&1
\end{array}\right)
\]
which at $c=-b=1$ reproduces eqn~\eqref{eq:ro2}.
Element $\theta$ is unipotent when $bc=1$, and we have
\[
j(x)=(\theta
g)_+=\left(\begin{array}{cc}-bx_1x_2^{-1}&b^2x_2^{-1}\\0&1
\end{array}\right)
\]
which at $b=1$ reproduces Example~\ref{ex:2}.
\end{example}
\begin{example}\label{ex:sl3}
Let $G=SL(3,\COMPLEXS)$ and $s_1,s_2,s_3$ be as in Example~\ref{ex:3}.
A symmetric factorization is defined
by the subgroups
\begin{equation}\label{eq:g+}
G_+=\left\{\left(\begin{array}{ccc} x_0^{s_1}& x_0^{s_1}x_1&x_0^{s_1}x_2\\
                               0&x_0^{s_2}& x_0^{s_2}x_3\\
                               0&0& x_0^{s_3}
                               \end{array}\right)|\
(x_0,x_1,x_2,x_3)\in\COMPLEXS^*\times\COMPLEXS^3\right\}
\end{equation}
\begin{equation}\label{eq:g-}
G_-=\left\{\left(\begin{array}{ccc} y_0^{s_3}& 0&0\\
                               y_0^{s_2}y_3&y_0^{s_2}& 0\\
                               y_0^{s_1}y_2&y_0^{s_1}y_1& y_0^{s_1}
                               \end{array}\right) |\
(y_0,y_1,y_2,y_3)\in \COMPLEXS^*\times\COMPLEXS^3\right\} .
\end{equation}
On a Zariski open subset of $G$ given by $g_{33}\bar g_{11}\ne0,\
\bar g\equiv g^{-1},$ there exists a unique factorization
\[
g=\|g_{ij}\|_{i,j=1}^3=g_+g_-^{-1}
\]
with
\begin{equation}\label{coord}
(g_+)^{-1}\leftrightarrow
 \left(\left(\bar g_{11}^{s_1}g_{33}^{s_3}\right)^{1/(s_2s_{31})},
\frac{\bar g_{12}}{\bar g_{11}},\frac{\bar g_{13}}{\bar g_{11}},
-\frac{g_{23}}{g_{33}}\right)
\end{equation}
and
\[
(g_-)^{-1}\leftrightarrow
 \left(\left(\bar g_{11}^{s_3}g_{33}^{s_1}\right)^{1/(s_2s_{31})},
\frac{g_{32}}{g_{33}},\frac{g_{31}}{g_{33}}, -\frac{\bar
g_{21}}{\bar g_{11}}\right) .
\]
Here four-component vectors parametrize matrices $(g_\pm)^{-1}$
according to eqns~\eqref{eq:g+} and \eqref{eq:g-}. For any
$a,b\in\COMPLEXS^*$ the element
\[
\theta=\left(\begin{array}{ccc}0&0&b\\0&-a/b&0\\
1/a
&0&0\end{array}\right)
\]
is conjugating. It is unipotent if $a=b$.

Set $x_4= x_1x_3-x_2$, where $(x_0,x_1,x_2,x_3)$ are coordinates on
$G_+$ (as in ~\eqref{coord}), then the inverse mapping $i(x)$ has
the form
\[
i(x)=(x_0^{-1},-x_0^{s_{12}}x_1,x_0^{s_{13}}x_4, -x_0^{s_{23}}x_3)
.
\]
The birational correspondence $\rho(x)=(\theta x)_+^{-1}$ is
given by
\[
\rho(x)=\left(
\left((x_2/a)^{s_3}(x_4/b)^{s_1}\right)^{1/(s_2s_{31})},
a^{-1}b^2x_0^{s_{32}}x_1x_4^{-1},abx_0^{s_{31}}x_4^{-1},
a^2b^{-1}x_0^{s_{21}}x_3x_2^{-1} \right) ,
\]
which at $a=b=1$ reproduces eqn~\eqref{eq:rho-sl3}. When $a=b$,
the correspondence $j(x)=(\theta x)_+$ is of order two and can be
written as a composition $j=iki=kik$, where
\[
k(x)= \left(x_0\left(
b^{s_2}x_2^{s_1}x_4^{s_3}\right)^{1/(s_2s_{31})},
-bx_1x_2^{-1},b^2 x_2^{-1},-b x_3x_4^{-1}\right) .
\]
\end{example}
Thus, according to Proposition \ref{prop:rho}, on the group $G_+$
we have the group multiplication $\cdot$ and  the binary operation
given by the  birational correspondence $*$ defined by the
factorizable structure on $G_+$. Therefore, for each symmetrically
factorizable algebraic Lie group we have a set-theoretical
solution to the pentagon equation.

Finally, one should notice that by varying $\theta$, one can have
a family of set-theoretical solutions to the pentagon equation.
Let $h\in G$ be an element such that $hG_-= G_- h$. Then, $\theta
h$ is again a conjugating element and, therefore, defines a $*$
operation on $G_+$. This gives a family of solutions to the
pentagon equation defined by a symmetrically factorizable Lie
group. We have such families in Examples~\ref{ex:gl2} and
\ref{ex:sl3}.

We have chosen  the algebro-geometrical setting, because it is the
case when, in a forthcoming paper we will use our results for
studying the moduli space of flat $G$-connections. Similar
analysis can be done in other topological settings. We will not do
it here.

\end{document}